\documentclass[11pt]{article} \topmargin=-0.5cm
\usepackage{amssymb}

   \csname @addtoreset\endcsname{equation}{section}
    \csname @addtoreset\endcsname{figure}{section}
    \csname @addtoreset\endcsname{table}{section}

\newcounter{thanksnum}
\def\thanksnumber#1
{\setcounter{thanksnum}{\value{footnote}}\setcounter{footnote}{#1}%
                     \addtocounter{footnote}{-1}\footnotemark
                     \setcounter{footnote}{\value{thanksnum}}}
\def\newtheoremz#1{\@ifnextchar[{\@othmz{#1}}{\@nthmz{#1}}}

\def\@nthmz#1#2{%
\@ifnextchar[{\@xnthmz{#1}{#2}}{\@ynthmz{#1}{#2}}}

\def\@xnthmz#1#2[#3]{\expandafter\@ifdefinable\csname #1\endcsname
{\@definecounter{#1}\@addtoreset{#1}{#3}%
\expandafter\xdef\csname the#1\endcsname{\expandafter\noexpand
  \csname the#3\endcsname \@thmcountersepz \@thmcounterz{#1}}%
\global\@namedef{#1}{\@thmz{#1}{#2}}\global\@namedef{end#1}{\@endtheoremz}}}

\def\@ynthmz#1#2{\expandafter\@ifdefinable\csname #1\endcsname
{\@definecounter{#1}%
\expandafter\xdef\csname the#1\endcsname{\@thmcounterz{#1}}%
\global\@namedef{#1}{\@thm{#1}{#2}}\global\@namedef{end#1}{\@endtheoremz}}}

\def\@othmz#1[#2]#3{\expandafter\@ifdefinable\csname #1\endcsname
  {\global\@namedef{the#1}{\@nameuse{the#2}}%
\global\@namedef{#1}{\@thmz{#2}{#3}}%
\global\@namedef{end#1}{\@endtheoremz}}}

\def\@thmz#1#2{\refstepcounter
    {#1}\@ifnextchar[{\@ythmz{#1}{#2}}{\@xthmz{#1}{#2}}}

\def\@xthmz#1#2{\@begintheoremz{#2}{\csname the#1\endcsname}\ignorespaces}
\def\@ythmz#1#2[#3]{\@opargbegintheoremz{#2}{\csname
       the#1\endcsname}{#3}\ignorespaces}

\def\@thmcounterz#1{\noexpand\arabic{#1}}
\def\@thmcountersepz{.}
\def\@begintheoremz#1#2{ \trivlist \item[\hskip \labelsep{\bf #1\ #2}]}
\def\@opargbegintheoremz#1#2#3{ \trivlist
      \item[\hskip \labelsep{\bf #1\ #2\ (#3)}]}
\def\@endtheoremz{\endtrivlist}

\newtheorem{theorem}{Theorem}[section]
\newtheorem{lemma}{Lemma}[section]

\newtheorem{corollary}{Corollary}[section]
\newtheorem{condition}{Condition}[section]

\def\e{\varepsilon}

\def\defi{\stackrel{{\scriptscriptstyle \Delta}}{=}}

\def\d{\delta}

\def\w{\widehat}
\def\Ind{{\mathbb{I}}}

\def\esssup{\mathop{\rm ess\, sup}}
\def\ulim{\mathop{\overline{\rm lim}}\,}

\def\R{{\bf R}}

\def\b{\beta}

\def\g{\gamma}

\def\ww{\widetilde}

\def\t{\theta}
\def\oo{\bar}

\def\p{\partial}

\def\A{{\cal A}}

\newcommand{\be}{\begin{equation}}
\newcommand{\ee}{\end{equation}}
\newcommand{\bd}{\begin{displaymath}}
\newcommand{\ed}{\end{displaymath}}
\newcommand{\ba}{\begin{array}{ll}}
\newcommand{\ea}{\end{array}}
\newcommand{\baa}{\begin{eqnarray}}
\newcommand{\eaa}{\end{eqnarray}}
\newcommand{\baaa}{\begin{eqnarray*}}
\newcommand{\eaaa}{\end{eqnarray*}}   \font\sm=cmr10

\def\ww{\tilde}

\def\CC{{\cal C}}

\def\FF{F}
\title{Universal  estimates for
parabolic equations and applications for non-linear and non-local
problems }
\author{
Nikolai Dokuchaev
\\ {\sm Department of Mathematics, Trent University, Ontario,
Canada}}
\begin{document}
\maketitle
\begin{abstract}
We obtain some "universal" estimates for $L_2$-norm of the
solution of a parabolic equation via a weighted version of
$H^{-1}$-norm of the free term. More precisely, we found the limit
upper estimate  that can be achieved by transformation of the
equation by adding a constant to the zero order coefficient. The
inverse matrix of the higher order coefficients of the parabolic
equation is included into the weight for the $H^{-1}$-norm. The
constant in the estimate obtained is independent from the choice
of the dimension, domain, and the coefficients of the parabolic
equation, it is why it can be called an universal estimate. As an
example of applications, we found an asymptotic upper estimate for
the norm of the solution at initial time. As an another example,
we established existence and regularity for non-linear and
non-local problems.
\end{abstract}
\section{Introduction}
We study prior estimates for first boundary value problems for
parabolic equations. The classical results for these equation give
upper estimate for the $L_2$-type Sobolev norm of the solution via
a $H^{-1}$-norm of the nonhomogeniuos term, where $H^{-1}$ is the
space being dual to the space $\stackrel{\scriptscriptstyle
0}{W_2^1}(D)$ (see, e.g., the first energy inequality  in
Ladyzhenskaia (1985)). We suggest a modification of this estimate.
\par
We found the limit minimal upper estimate that can be achieved by
varying the zero order coefficient of the original equation by
adding a constant. In other words, we study  the case when the
original equation is transformed into a new one such that the
original solution $u(x,t)$ is to be replaced by $u(x,t)e^{-Kt}$;
the value of $K$ is being varied (Theorem \ref{ThM} and Lemma
\ref{lemma<<1}). The constant in the estimate is the same  for all
possible choices of the dimension, domain, time horizon, and the
coefficients of the parabolic equations. It is why it can be
called a universal estimate. These results represent an important
development of the extension of the results from Dokuchaev (2008),
where an "universal" estimate was obtained for the gradient via
$L_2$-norm of the nonhomogeniuos term. In contrast, the present
paper gives the estimate of the $L_2$-norm via a $H^{-1}$-type
norm of the nonhomogeniuos term, i.e., via a weaker norm. It is
shown that the estimate obtained is sharp (Theorem
\ref{ThSharpM}).
\par
 As an example of applications, we obtained  a sharp asymptotic upper
estimate for the solution at initial time (Theorems \ref{ThAs} and
\ref{ThSharp}). The constant in this eastimate is again the same
for all possible equations. As an another example of applications,
we suggest a new approach for establishing of existing and
regularity for non-linear and non-local parabolic equations
(Theorem \ref{Thn}). We found an explicit sufficient conditions
for existence and regularity (Conditions (\ref{LipY1}) or
(\ref{LipY2})). These conditions are easy to verify, and they
cover a wide class of non-linear and non-local parabolic
equations.
\section{Definitions}\label{SecD}
\subsection*{Spaces and classes of functions.} 
 We denote by $|\cdot|$  the Euclidean norm in $\R^k$ and the Frobenius norm in $\R^{k\times m}$,
 and we denote by $\bar G$ denote the closure of a region $G\subset\R^k$.
\par
\par
We denote by $\|\cdot\|_{ X}$ the norm in a linear normed space
$X$, and
 $(\cdot, \cdot )_{ X}$ denote  the scalar product in  a Hilbert space $
X$. For a Banach space $X$, we denote by $C([a,b],X)$ the Banach
space of continuous functions $x:[a,b]\to X$.
\par
 Let $G\subset \R^k$ be an open
domain, then ${W_q^m}(G)$ denote  the Sobolev  space of functions
that belong $L_q(G)$ together with the distributional derivatives
up to the $m$th order, $q\ge 1$.
\par
We are given an open domain $D\subseteq\R^n$ such that either
$D=\R^n$ or $D$ is bounded  with $C^2$-smooth boundary $\p D$.
\par
Let $T>0$ be given, and let $Q\defi D\times (0,T)$.
\par Let $H^0\defi L_2(D)$,
and let $H^1\defi \stackrel{\scriptscriptstyle 0}{W_2^1}(D)$ be the
closure in the ${W}_2^1(D)$-norm of the set of all smooth functions
$u:D\to\R$ such that  $u|_{\p D}\equiv 0$. The spaces $H^k$ are
Hilbert spaces, and $H^1$ is a closed subspace of $W_2^1(D)$.  Let
$H^2=W^2_2(D)\cap H^1$ be the space equipped with the norm of
$W_2^2(D)$. The spaces $H^k$ are Hilbert spaces, and $H^k$ is a
closed subspace of $W_2^k(D)$,  $k=1,2$.
\par
 Let $H^{-1}$ be the dual space to $H^{1}$, with the
norm $\| \,\cdot\,\| _{H^{-1}}$ such that if $u \in H^{0}$ then
$\| u\|_{ H^{-1}}$ is the supremum of $(u,w)_{H^0}$ over all $w
\in H^1$ such that $\|w\|_{H^1} \le 1 $. $H^{-1}$ is a Hilbert
space.
\par We will write $(u,w)_{H^0}$ for $u\in H^{-1}$
and $w\in H^1$, meaning the obvious extension of the bilinear form
from $u\in H^{0}$ and $w\in H^1$.
\par
We denote by $\oo\ell _{1}$ the Lebesgue measure in $\R$, and we
denote by $ \oo{{\cal B}}_{1}$ the $\sigma$-algebra of Lebesgue sets
in $\R^1$.\
\par
 For $k=-1,0,1,2$, we  introduce the spaces
 \baaa
 X^{k}\defi L^{2}\bigl([ 0,T ], \oo{{\cal B}}_{1},\oo\ell_{1};  H^{k}\bigr),
 \quad \CC^{k}\defi C\left([0,T]; H^k\right). \eaaa
\par
We introduce the spaces $$ Y^{k}\defi X^{k}\!\cap \CC^{k-1}, \quad
k=0,1,2, $$ with the norm $ \| u\| _{Y^k}
\defi \| u\| _{{X}^k} +\| u\| _{\CC^{k-1}}. $
\par
\par
 We use the notations $$ \nabla u\defi \Bigl(\frac{\p u}{\p
x_1},\frac{\p u}{\p x_2},\ldots,\frac{\p u}{\p
x_n}\Bigr)^\top,\quad \nabla\cdot U=\sum_{i=1}^n\frac{\p U_i}{\p
x_i} $$ for functions $u:\R^n\to\R$ and
$U=(U_1,\ldots,U_n)^{\top}:\R^n\to \R^n$. In addition, we  use the
notation $$ (U,V)_{H^0}=\sum_{i=1}^n(U_i,V_i)_{H^0},\quad
\|U\|_{H^0}= (U,U)_{H^0}^{1/2} $$ for functions $U,V: D\to\R^n$,
where $U=(U_1,\ldots,U_n)$ and $V=(V_1,\ldots,V_n)$.
\subsubsection{The boundary value problem}
We consider the following problem \be \label{4.1} \ba \frac{\p
u}{\p t} =\A u+ \varphi,\quad \quad t\in( 0,T),\\ u(x,0)=0,\quad
u(x,t)|_{x\in \p D}=0 . \ea
 \ee
 Here
 $u=u(x,t)$,
 $(x,t)\in Q$,  and
  \be\label{A}\A y\defi  \sum_{i=1}^n\frac{\p }{\p x_i}
\sum_{j=1}^n \Bigl(b_{ij}(x,t)\frac{\p y}{\p x_j}(x)\Bigr)
+\sum_{i=1}^n f_{i}(x,t)\frac{\p y}{\p x_i }(x)
+\,\lambda(x,t)y(x), \ee where $b(x,t):
\R^n\times[0,T]\to\R^{n\times n}$, $f(x,t):
\R^n\times[0,T]\to\R^n$, and $\lambda(x,t): \R^n\times[0,T]\to\R$,
are bounded measurable functions, and $b_{ij}, f_i, x_i$ are the
components of $b$,$f$, and $x$. The matrix $b=b^\top$ is
symmetric.
\par
 To proceed further, we assume that Conditions
\ref{cond3.1.A}-\ref{cond3.1.B} remain in force throughout this
paper.
 \begin{condition} \label{cond3.1.A}   There exists a constant $\d>0$ such that
\be
 \label{Main1} \xi^\top  b
(x,t)\,\xi \ge \d|\xi|^2 \quad\forall\, \xi\in \R^n,\ (x,t)\in Q.
\ee
\end{condition}
\par
Inequality (\ref{Main1})  means that equation (\ref{4.1}) is
coercive.
\begin{condition}\label{cond3.1.B} The
functions $b(x,t):\R^n \times \R\to \R^{n \times n}$, $f(x,t):\R^n
\times \R\to \R^n$, $\lambda (x,t):\R^n \times \R\to \R$, are
measurable, and \baaa \esssup_{(x,t)\in Q} \biggl[| b(x,t)|+|
b(x,t)^{-1}|+ |f(x,t)|+ | \lambda(x,t)|\biggl]< +\infty. \eaaa
\end{condition}
\par
We introduce the sets of parameters   \baaa
 &&\mu
\defi (T,\, n,\,\, D,\,\,  \, \delta,\,\,\,\,
v, f,\lambda),
 \\
 &&{\cal P}={\cal P}(\mu)
\defi \biggl(T,\, n,\,\, D,\,\,   \delta,\,\,\,\,
\esssup_{(x,t)\in Q}\Bigl[| b(x,t)|+ |f(x,t)|+ |
\lambda(x,t)|\Bigr] \biggr). \eaaa We consider all possible $\mu$
such that the conditions imposed above are satisfied.
 \section{Special estimates for the solution}
We assume that $\varphi\in X^{-1}$. This means that there exist
functions $\FF=(\FF_1,...,\FF_n):Q\to\R^n$ and $\FF_0:Q\to\R$ such
that $\FF_k\in X^0=L_2(Q)$, $k=0,1,...,n$, and
\be\varphi(x,t)=\nabla\cdot \FF(x,t)+\FF_0(x,t).\label{hU}\ee In
other words, $\varphi(x,t)=\sum_{k=1}^n\frac{\p \FF_k}{\p
x_k}(x,t)+\FF_0(x,t)$.
\par
 The classical solvability results for
the parabolic equations give that there exists a unique solution
$u\in Y^1$ of problem (\ref{4.1}) for any $\varphi\in X^{-1}$.
 In addition, it follows from the first
energy inequality  (or the first fundamental inequality) that, for
any $K\in \R$ and $M\ge 0$, there exist constants $\ww
C_i(K,M,{\cal P})>0$, $i=0,1$, such that \baa
&&e^{-2Kt}\|u(\cdot,t)\|^2_{ H^0}+M\int_0^t
e^{-2Ks}\|u(\cdot,s)\|^2_{ H^0}ds \nonumber\\ &&\hspace{1cm}\le\ww
C_1(K,M,{\cal P})\int_0^t
e^{-2Ks}(\FF(\cdot,s),b(\cdot,s)^{-1}\FF(\cdot,s))_{H^0}ds
\nonumber\\ &&\hspace{2cm}+\ww C_0(K,M,{\cal P})\int_0^t
e^{-2Ks}\|\FF_0(\cdot,s)\|_{H^0}^2ds\qquad \forall \varphi\in
X^{-1},\ t\in (0,T], \label{esti}\eaa where $F_i\in X^0$ are such
that (\ref{hU}) holds. (See, e.g., estimate (3.14) from
Ladyzhenskaia (1985), Chapter III, \S 3). We have used here the
following obvious estimate
$$\sum_{k=1}^n\|\FF_k(\cdot,s)\|_{H^0}^2ds \le
c(\FF(\cdot,s),b(\cdot,s)^{-1}\FF(\cdot,s))_{H^0}, $$ where
$c=c({\cal P})>0$ is a constant.
\par Let
$C_i(K,M,{\cal P})\defi\inf \ww C_i(K,M,{\cal P})$,  where the
infimum is taken over  all  $\ww C_i(K,M,{\cal P})$ such that
(\ref{esti}) holds.
\begin{theorem}\label{ThM}  $$\sup_{\mu,\,M\ge 0} \inf_{K\ge 0} C_1(K,M,{\cal
P}(\mu))\le \frac{1}{2},\quad \sup_{\mu,\,M\ge 0} \inf_{K\ge 0}
C_0(K,M,{\cal P}(\mu))= 0. $$
\end{theorem}
\begin{corollary}\label{corr1}
For any $\mu$ and any $M>0$, $\e>0$, there exists $ K=
K(\e,M,{\cal P(\mu)})\ge 0$ such that \baa
&&\sup_{s\in[0,t]}e^{-2Ks}\|u(\cdot,s)\|^2_{ H^0}+M\int_0^t
e^{-2Ks}\|u(\cdot,s)\|^2_{ H^0}ds\nonumber\\&&\hspace{5mm}\le
\Bigl(\frac{1}{2}+\e\Bigr)\sum_{k=0}^n\int_0^t(\FF(\cdot,s),b(\cdot,s)^{-1}\FF(\cdot,s))_{H^0}ds
+\e \int_0^t\|\FF_0(\cdot,s)\|_{H^0}^2ds \label{sup}\\\ \quad
&&\hspace{6cm}\forall t\in [0,T],\ \varphi\in X^{-1},\nonumber
\eaa
 where $u$ is the solution
 of problem  (\ref{4.1}), and  where $F_i\in X^0$ are such
that (\ref{hU}) holds.
\end{corollary}
\section{The case of non-linear and non-local equations}
Let us consider the following mapping ${\cal N}(v):Y^1\to X^{-1}$
such that
  \baa{\cal N}(v) \defi  \sum_{i=1}^n\frac{\p }{\p x_i}
\sum_{j=1}^n \Bigl(\w b_{ij}(v(\cdot),x,t)\frac{\p v}{\p
x_j}(x,t)\Bigr) +\sum_{i=1}^n\w f_{i}(v(\cdot),x,t)\frac{\p v}{\p
x_i }(x,t)\nonumber\\ +\,\w\lambda(v(\cdot),x,t)v(x,t)
+\w\varphi(v(\cdot),x,t),\hphantom{} \label{An}\eaa where $\w
b(v(\cdot),x,t): Y^1\times Q\to\R^{n\times n}$, $\w
f(v(\cdot),x,t): Y^1\times Q\to\R^n$, $\w\lambda(v(\cdot),x,t):
Y^1\times Q\to\R$, are bounded functions, and $b_{ij}, f_i, x_i$
are the components of $b$,$f$, and $x$. The function
$\w\varphi(v(\cdot),x,t)$ is defined on $Y^1\times Q\to\R$ and
belongs to $X^{-1}$ for any given $v(\cdot)\in Y^1$.
 The
matrix $\w b=\w b^\top$ is symmetric.
\begin{corollary}
\label{corrn} Let $u\in Y^1$ be a solution of the problem \be
\label{4.1n} \ba \frac{\p u}{\p t}
 ={\cal N}(u),\quad \quad t\in( 0,T),\\ u(x,0)=0,\quad u(x,t)|_{x\in
\p D}=0. \ea
 \ee
 such that Conditions \ref{cond3.1.A}-\ref{cond3.1.B} are satisfied for
 $b(x,t)=\w b(u(\cdot),x,t)$, $f(x,t)=\w f(u(\cdot),x,t)$, and  $\lambda(x,t)=\w
 \lambda(u(\cdot),x,t)$, and such that  $\varphi(x,t)\defi
\w\varphi(u(\cdot),x,t)$ belongs to $X^{-1}$ and is such that
(\ref{hU}) holds for $F_i\in X^0$. Then, for any $M>0$ and $\e>0$,
there exists $ K= K(\e,M,{\cal P})\ge 0$ such that (\ref{sup})
holds.
\end{corollary}
\par
Note that  the parabolic equation in (\ref{4.1n}) is non-linear
and non-local.
\par
Corollary \ref{corrn} does not establish existence. Some existence
results for non-local and non-linear problems are  given below.
\section{Applications: asymptotic estimate at initial time}
Let \baa X^0_c\defi\Bigl\{\varphi\in X^{0}:\quad \lim_{t\to
0+}\frac{1}{t}\int_0^t\|\varphi(\cdot,s)\|_{H^0}^2ds
=\|\varphi(\cdot,0)\|_{H^0}^2\Bigr\}.\label{Xc0}\eaa  Note that
the condition that $\varphi\in X^{0}_c$ is not restrictive for
$\varphi\in X^{0}$; for instance, it holds if $s=0$ is a Lebesgue
point for $\|\varphi(\cdot,s)\|_{H^0}^2$.
\begin{theorem}\label{ThAs0}
Let $\varphi\in X_c^{0}$. Then, for any admissible $\mu$, the
solution $u$
 of problem  (\ref{4.1})
 is such that \baa\ulim_{t\to 0+}\sup_{\varphi\in X_c^{-1}}\frac{1}{t}\frac{\|u(\cdot,t)\|^2_{
H^0}}{\|\varphi(\cdot,0)\|_{H^0}^2}=0, \label{as10}\eaa where $u$
is the solution
 of problem  (\ref{4.1}) for the corresponding $\varphi$.
\end{theorem}
\par
Further, let \baa X^{-1}_c\defi\Bigl\{\varphi\in X^{-1}:\hbox{
there exists a set $\{\FF_k\}_{k=1}^n\subset X^0$ such that
(\ref{hU}) holds with $F_0\equiv 0$,} \nonumber\\\hbox{and }
\lim_{t\to
0+}\frac{1}{t}\int_0^t(\FF(\cdot,s),b(\cdot,s)^{-1}\FF(\cdot,s))_{H^0}ds
=(\FF(\cdot,0),b(\cdot,0)^{-1}\FF(\cdot,0))_{H^0}\Bigr\}.\hphantom{x}\label{Xc}\eaa
Here $F=(F_1,...,F_n)$. Again, the limit condition in (\ref{Xc})
is not restrictive; for instance, it holds if $s=0$ is a Lebesgue
point for $(\FF(\cdot,s),b(\cdot,s)^{-1}\FF(\cdot,s))_{H^0}$.
Clearly, the set $X_c^{-1}$ is non-empty if some mild conditions
of regularity in $t$ are satisfied for   $b(x,t)^{-1}$.
\begin{theorem}\label{ThAs}
Let $\varphi\in X_c^{-1}$ and let $\FF_k\in X^0$, $k=0,1,...,n$,
be the corresponding functions presented in (\ref{Xc}) with
$F_0=0$. Then, for any admissible $\mu$,  \baa\ulim_{t\to
0+}\sup_{\varphi\in X_c^{-1}}\frac{1}{t}\frac{\|u(\cdot,t)\|^2_{
H^0}}{(\FF(\cdot,0),b(\cdot,0)^{-1}\FF(\cdot,0))_{H^0}}\le
\frac{1}{2}, \label{as1}\eaa where $u$ is the solution
 of problem  (\ref{4.1}) for the corresponding $\varphi$, $F=(F_1,....,F_n)$.
\end{theorem}
\begin{corollary}\label{corrAs}
Let $u\in Y^1$ be a solution
 of problem  (\ref{4.1n}) such that the assumptions of Corollary \ref{corrn}
 are satisfied.
Let $\varphi(x,t)=\w\varphi(u(\cdot),x,t)$ be such that
$\varphi=\FF_0+\ww\varphi$, where
  $\FF_0\in X^0_c$ and
$\ww\varphi\in X_c^{-1}$, and let $\FF_i\in X^0$ be the
corresponding functions presented in (\ref{hU})  such that the
limit conditions from (\ref{Xc}) are satisfied. Then
 \baa\ulim_{t\to 0+}\frac{1}{t}\|u(\cdot,t)\|^2_{
H^0}\le
\frac{1}{2}(\FF(\cdot,0),b(\cdot,0)^{-1}\FF(\cdot,0))_{H^0},
\label{as101}\eaa where $\FF=(\FF_1,...,\FF_n)$.
\end{corollary}
Note that $\FF_0$ is not being presented in the last estimate.
 \section{Applications: existence for non-linear and non-local
 equations}The universal estimates from
 Theorem  \ref{ThM} can be also applied to analysis of non-linear and non-local parabolic
 equations. These equations have many applications, and they were intensively
 studied (see. e.g., Ammann, (2005),
 Ladyzenskaya {\it et al} (1967),
Zheng (2004), and references there). Theorem \ref{ThM} gives a new
way to establish  conditions of
  solvability of these equations. This approach
  covers many cases when  the solutions and the gradient are
 included into the non-local and non-linear term.
 \par
Let $B(u(\cdot)):X^0\to X^{-1}$ be a mapping that describes
non-linear and non-local term in the equation.
\par
 Let
us consider the following boundary value problem in $Q$: \be
\label{D} \ba \frac{\p u}{\p t} =\A u+ B(u)+\varphi,\quad \quad
t\in( 0,T),\\ u(x,0)=0,\quad u(x,t)|_{x\in \p D}=0 . \ea
 \ee
 Here $\A$ is the linear operator defined above.
 For $K>0$, introduce the mappings
\be\label{BK} B_K(u)\defi e^{-Kt}B(\oo u_K),\quad
\hbox{where}\quad \oo u_K(x,t)\defi e^{Kt}u(x,t). \ee
 \begin{theorem}\label{Thn}  Assume that $B(u)$ maps $X^0$ into $X^{-1}$.
 Moreover, assume that there exist constants $K_*>0$ and $C_*>0$
such that
 \baa
 \|B_K(u_1)-B_K(u_2)\|_{X^{-1}}\le
 C_*\|u_1-u_2\|_{X^{0}}
 \quad \forall u_1,u_2\in X^{0},\ K>K_*.
 \label{LipY1}
 \eaa
 Then there exists a unique solution $u\in Y^1$ of problem
(\ref{D}) for any $\varphi\in X^{-1}$. \index{and $\|u\|_{Y^1}\le
c\|\varphi\|_{X^{-1}}$, where $c$ is a constant that depends only on
${\cal P}$ and $C_*$.}
\end{theorem}
  \begin{theorem}\label{ThnY2}  Assume that $B(u)$ maps $X^1$ into $X^{0}$
  and that there exist constants
$K_*>0$ and $C_*>0$ such that
 \baa
 \|B_K(u_1)-B_K(u_2)\|_{X^{0}}\le
 C_*\|u_1-u_2\|_{X^{1}}
 \quad \forall u_1,u_2\in X^{1},\ K>K_*.
 \label{LipY2}
 \eaa
 Then there exists a unique solution $u\in Y^2$ of problem
(\ref{D}) for any $\varphi\in X^{0}$.
\end{theorem}
\subsubsection{Examples of admissible $B$}
 Some examples covered by Theorem \ref{Thn} are listed below.
 \begin{theorem}\label{ThEx}
 The assumptions of Theorem \ref{Thn}  hold for the following
 mappings $B(u)$:
 \begin{itemize}\item[(i)] A local non-linearity:
 \be\label{Bsin} B(u)=\b (u(x,t),x,t), \ee
 where $\b: \R\times Q\to\R$ is a
 measurable  function such that $\b(0,\cdot)\in L_2(Q)$ and that there exists a constant $C_L>0$ such
 that
\baa
 |\b(z_1,x,t)-\b(z_2,x,t)|\le
 C_L|z_1-z_2|
 \quad \forall z_1,z_2\in\R,\ x,t. \label{Lipii}\eaa
 \item[(ii)]  A distributional non-linearity:
 \be\label{Bsin'} B(u)\defi \nabla\cdot \b (u(x,t),x,t), \ee
 where $\b: \R\times Q\to\R^n$ is a
 measurable function such that  $\b(0,\cdot)\in L_2(Q)$  and (\ref{Lipii}) holds.
\item[(iii)] A non-local non-linearity (integral nonlinearity):
 \baaa\label{B30} (B(u))(x,t)=\int_D\b(u(y,t),x,t,y)dy, \eaaa
 where $\b: \R\times Q\times D\to\R$ is a
 measurable function such that
 $\int_D\b(0,x,t,y)dy\in L_2(Q)$ as a function of
 $(x,t)$,
 and there exists a constant $C_L>0$ such
 that
 \baa |\b(z_1,x,t,y)-\b(z_2,x,t,y)|\le
 C_L|z_1-z_2|\quad \forall z_1,z_2\in\R,\ x,t,y.
 \label{Lip220}\eaa We assume here that $D$ is a bounded domain.
 \item[(iv)] A non-local in space distributional non-linearity:
 \baaa\label{B3'0} (B(u))(x,t)= \nabla\cdot\int_D\b(u(y,t),x,t,y)dy, \eaaa
 where $\b: \R\times Q\times D\to\R^n$ is a
 measurable function such that $\int_D\b(0,\cdot,y)dy\in
 L_2(Q)$ as a function of $(x,t)$,
 and (\ref{Lip220}) holds. We assume here that $D$ is a bounded domain.
 \item[(v)] A non-local  in time and space non-linearity:
 \baaa\label{B3} (B(u))(x,t)=\int_0^tds\int_D\b(u(y,s),x,t,y,s)dy, \eaaa
 where $\b: \R\times Q^2\to\R$ is a
 measurable function such that
 $\int_0^tds\int_D\b(0,x,t,y,s)dy\in L_2(Q)$ as a function of
 $(x,t)$,
 and there exists a constant $C_L>0$ such
 that
 \baa |\b(z_1,x,t,y,s)-\b(z_2,x,t,y,s)|\le
 C_L|z_1-z_2|\quad \forall z_1,z_2\in\R,\ x,t,y,s.
 \label{Lip22}\eaa We assume here that $D$ is a bounded domain.
 \item[(vi)] A non-local  in time and space distributional non-linearity:
 \baaa\label{B3'} (B(u))(x,t)= \nabla\cdot\int_0^tds\int_D\b(u(y,s),x,t,y,s)dy,
 \eaaa
 where $\b: \R\times Q^2\to\R^n$ is a
 measurable function such that $\int_0^tds\int_D\b(0,\cdot,y,s)dy\in
 L_2(Q)$ as a function of $(x,t)$,
 and (\ref{Lip22}) holds. We assume here that $D$ is a bounded domain.
 \item[(vii)] Nonlinear delay parabolic equations:
  \be\label{Bd}  (B(u))(x,t)\defi\nabla\cdot\beta(u(x,\tau(t)),x,\tau(t))
  +\oo \beta(u(x,\tau(t)), x,\tau(t)).
  \ee
 Here
$\tau(\cdot):[0,T]\to\R$ is a given measurable function such that
$\tau(t)\in [0,t]$,  and that there exists $\t\in[0,T)$ such that
$\tau(t)=0$ for $t<\t$, the function $\tau(\cdot):[\t,T]\to \R$ is
non-decreasing and absolutely continuous, and
 $\esssup_{t\in[\t,T]}\left|\frac{d\tau}{dt}(t)\right|^{-1}<+\infty$.
The functions $\b:\R\times \R^n\times[0,T]\to\R^n$ and
$\oo\b:\R\times\R^n\times[0,T]\to\R$ are bounded and measurable.
In addition, we assume that the derivative $\frac{\p \b}{\p
x}(x,t)$ is bounded,
 $\b(0,\cdot)\in L_2(Q)$,  $\w\b(0,\cdot)\in L_2(Q)$, and there exists a constant $C_L>0$ such
 that
 \baa |\b(z_1,x,t)-\b(z_2,x,t)|+ |\w\b(z_1,x,t)-\w\b(z_2,x,t)|\le
 C_L|z_1-z_2|\quad \forall z_1,z_2\in\R,\ x,t.\hphantom{x}
  \label{LipD}\eaa
\item[(viii)] Non-local term for the backward Kolmogorov equations for a jump
diffusion process:
 \baaa\label{Bjump} (Bu)(x,t)\defi \int_{\R^n}
\Ind_{\{x+c(x,y,t)\in D\}}
 (u(x+c(x,y,t),t)-u(x,t)-c(x,y,t)^\top\nabla
 u(x,t))\rho(y,t)dy.\hphantom{x}
 \eaaa
Here $\rho(y,t):\R^n\times[0,T]\to \R$ is a function such that
$\rho(\cdot)\in L_{\infty}([0,T],\ell_1,\oo{\cal B}_1,L_1(\R^n))$.
The function $c(x,y,t): D\times \R^n\times[0,T]\to\R^n$ is
measurable,
 bounded, and  such that the derivative  $\frac{\p c}{\p x}(x,y,t)$
 is bounded, the derivative  $\frac{\p c}{\p z}(x,y,t)$ exists almost everywhere,
  and there exists a uniquely defined
 function $\psi:D\times \R^n\times [0,T]\to\R^n$ such that $z=x+c(x,y,t)$ for
 $y=\psi(x,z,t)$.
In addition, we assume that  $\esssup_{t\in [0,T]} \int_{D\times
D} |r(x,z,t)|^2dxdz<+\infty$, where
$r(x,z,t)\defi\rho(\psi(x,z,t),t)\frac{\p \psi}{\p z}(x,z,t).$
 \end{itemize}
\end{theorem}
\par Clearly, linear combinations of the
non-linear and non-local terms listed above are also covered, as
well as terms formed as compound mappings.
\section{On the sharpness of the estimates}
\begin{theorem}\label{ThSharpM} There exists a set of parameters $(n,D,
b(\cdot), f(\cdot),\lambda(\cdot))$  such that, for any $T>0$,
$M\ge 0$, \be \inf_{K\ge 0} C(K,M,{\cal P}(\mu))= \frac{1}{2}.
\label{prot1}\ee for $\mu=(T,n,D,b(\cdot),
f(\cdot),\lambda(\cdot))$.
\end{theorem}
\begin{theorem}\label{ThSharp} There exists  a set of parameters $(n,D,
b(\cdot), f(\cdot),\lambda(\cdot))$   such that  \baa\ulim_{t\to
0+}\sup_{\varphi\in X_c^{-1}}\frac{1}{t}\frac{\|u(\cdot,t)\|^2_{
H^0}}{(\FF(\cdot,0),b(\cdot,0)^{-1}\FF(\cdot,0))_{H^0}}=
\frac{1}{2}.
 \label{as2}\eaa where
$u$ is the solution
 of problem  (\ref{4.1}) for the corresponding $\varphi\in X_c^{-1}$, and where $F=(F_1,...,F_n)$ with
 $\FF_i\in X^0$ being the corresponding functions  presented in (\ref{Xc}).
 \end{theorem}
\section{Proofs}
\begin{lemma}\label{lemma<<1} For any
admissible $\mu$ and any $\e>0$, $M>0$, there exists $\ww K=\ww
K(\e,M,{\cal P(\mu)})\ge 0$ such that \baa
&&\|u(\cdot,t)\|_{H^0}^2+ M\int_0^t
\|u(\cdot,s)\|_{H^0}^2ds\nonumber\\&&
\le\Bigl(\frac{1}{2}+\e\Bigr)\int_0^t\Bigl((\FF(\cdot,s),b(\cdot,s)^{-1}F(\cdot,s))_{H^0}ds
+\e \int_0^t\|\FF_0(\cdot,s)\|_{H^0}^2ds\eaa for all $K\ge \ww
K(\e,M,{\cal P})$, $t\in(0,T]$, for all $\varphi\in X^{-1}$
represented as (\ref{hU}) with  $\FF_i\in X^0$. Here $u\in Y^1$ is
the solution of the boundary value problem
 \be\ba \frac{\p u}{\p
t}=\A u-Ku+\varphi,\quad t\in(0,T),\\ u(x,0) =0,\quad
u(x,t)|_{x\in \p D}=0. \label{oog3}\ea\ee
\end{lemma}
\par
Uniqueness and existence of solution  $u\in Y^1$ of problem
(\ref{oog3}) follows from the classical results (see, e.g.,
Ladyzhenskaia (1985), Chapter III). \par {\it Proof of Lemma
\ref{lemma<<1}.} Clearly, $\A u=\A_su+\A_ru$, where \baaa \A_s
u=\nabla\cdot (b\nabla u)=\sum_{i=1}^n\frac{\p }{\p x_i}
\sum_{j=1}^n \Bigl(b_{ij}\frac{\p u}{\p x_j}\Bigr),\quad \quad\A_r
u=\sum_{i=1}^n f_i\frac{\p u}{\p x_i}+ \lambda u.
 \label{vv}\eaaa
Assume that $\varphi(\cdot,t)$ is differentiable and has a compact
support inside $D$ for all $t$.
 We have that \baa &&\|u(\cdot,t)\|^2_{H^0}-\|u(\cdot,0)\|^2_{H^0} =
 (u(\cdot,t),u(\cdot,t))_{H^0}-(u(\cdot,0),u(\cdot,0))_{H^0}\nonumber\\&& = 2\int_0^t \left(
u,\frac{\p u}{\p s}\right)_{H^0} ds  =2\int_0^t \left(u,\A
u-Ku+\varphi \right)_{H^0} ds  \nonumber\\&& =2\int_0^t \left(
u,\nabla\cdot(b\nabla  u\bigr)\right)_{H^0} ds +2\int_0^t \left(
u,\A_r u\right)_{H^0} ds-2K\int_0^t \left(u,u\right)_{H^0}
ds\nonumber\\&&\hspace{8cm}+2\int_0^t \left(u,
\varphi\right)_{H^0} ds.\hphantom{xx} \label{R0} \eaa
 Let arbitrary $\e_0>0$ and
$\w\e_0>0$ be given.
 Let $v\defi\sqrt{b}$, i.e.,
$b=v^2$, $v=v^\top$. We have that
 \baa &&2\left( u,\varphi\right)_{H^0} =2\left(u, \nabla\cdot \FF\right)_{H^0}
 +2\left(u, \FF_0\right)_{H^0}=-2\left(v\nabla u,v^{-1}\FF)\right)_{H^0}+2\left(u, \FF_0)\right)_{H^0}\nonumber\\
  &&\le
\frac{2}{1+2\e_0}\left(v\nabla u,  v\nabla
u\right)^2_{H^0}+\left(\frac{1}{2}+\e_0\right)\left\|v^{-1}\FF\right\|^2_{H^0}+
\frac{1}{\w\e_0}\left\|u\right\|^2_{H^0}+\w\e_0\left\|\FF_0\right\|^2_{H^0}
\nonumber\\
 &&=\frac{2}{1+2\e_0}\left(\nabla u, b \nabla
u\right)^2_{H^0}+\left(\frac{1}{2}+\e_0\right)\left(F,b^{-1}\FF\right)_{H^0}+
\frac{1}{\w\e_0}\left\|u\right\|^2_{H^0}+\w\e_0\left\|\FF_0\right\|^2_{H^0}
, \label{R1}\eaa and \baa 2\left( u,\nabla\cdot(b\nabla
u)\right)_{H^0} =-2\left(\nabla u,b\nabla u
\right)_{H^0}.\label{R2}\eaa In addition, we have that in under
the integrals in (\ref{R0}),
 \baaa 2\left( u, \A_r u\right)_{H^0} \le
\e_1^{-1}\left\|
u\right\|^2_{H^0}+\e_1\left\|\A_ru\right\|^2_{H^0}\quad\forall\e_1>0.
\label{R3'}\eaaa
\par
By the first energy inequality, there exist constants
$c'_{*}=c'_{*}({\cal P})>0$ and $c_{*}=c_{*}({\cal P})>0$ such
that \be\int_0^t\left\|u(\cdot,s)\right\|^2_{H^1}ds\le
c'_{*}\sum_{k=0}^n\int_0^t\left\|\FF_k (\cdot,s)\right\|^2_{H^0}ds
\le
c_{*}\int_0^t \left(F,b^{-1}\FF\right)_{H^0}ds. \label{e*} \ee
(See, e.g. inequality (3.14) from Ladyzhenskaia (1985), Chapter
III). Moreover, this constant $c_*$ can be taken the same for all
$t\in [0,T]$ and all $K>0$.
 Further, there exists a constant
$c_{1}=c_{1}({\cal P})>0$ such that  \baaa 2\left( u, \A_r
u\right)_{H^0} \le \e_1^{-1}\left\|
u\right\|^2_{H^0}+c_1\e_1\left\|u\right\|^2_{H^1}. \eaaa It
follows that \baa 2\int_0^t\left( u, \A_r u\right)_{H^0}ds \le
\e_1^{-1}\int_0^t\left\|
u\right\|^2_{H^0}ds+\e_0\int_0^t\left(F,b^{-1}\FF\right)_{H^0}ds,
\label{R3} \eaa
 if $\e_1>0$ is taken such that $c_1c_{*}\e_1=\e_0$.
\par
 By (\ref{R0})-(\ref{R3}), it follows that \baa
&& \|u(\cdot,t)\|^2_{H^0}+M\int_0^t \|u(\cdot,s)\|^2_{ H^0}ds
\nonumber\\&&\le \Bigl[\frac{2}{1+2\e_0}-2\Bigr]\int_0^t
\left(\nabla u,b\nabla u\right)_{H^0} ds
+[\e_1^{-1}+\w\e_0^{-1}+M-2K]\int_0^t \left\| u\right\|^2_{H^0} ds
\nonumber
\\&&
+\Bigl(\frac{1}{2}+2\e_0\Bigr)\int_{0}^t\left(F,b^{-1}\FF\right)_{H^0}ds
+(\e_0+\w\e_0)\int_0^t\left\|\FF_0 (\cdot,s)\right\|^2_{H^0}ds
\nonumber
\\&&\le
\Bigl(\frac{1}{2}+2\e_0\Bigr)\int_0^t\left(F,b^{-1}\FF\right)_{H^0}ds
+(\e_0+\w\e_0)\int_0^t\left\|\FF_0
(\cdot,s)\right\|^2_{H^0}ds,\hphantom{}\label{R4'} \eaa if
$2K>\e_1^{-1}+c_v'+M$. Then the proof of Lemma \ref{lemma<<1}
follows. $\Box$
 \par
{\it Proof of Theorem \ref{ThM}}. Clearly,
$u(x,t)=e^{Kt}u_K(x,t)$, where $u$ is the solution
 of problem  (\ref{4.1})  and $u_K$ is the solution of (\ref{oog3}) for
the  nonhomogeneous  term $e^{-Kt}\varphi(x,t)$. Therefore,
Theorem \ref{ThM} follows immediately from Lemma \ref{lemma<<1}.
$\Box$
 \par  Corollary \ref{corr1}
 follows immediately from Theorem \ref{ThM}.
 \par
 Corollary \ref{corrn}
 follows immediately from Corollary \ref{corr1}.
 \par {\it Proof of Theorem \ref{ThAs0}.}
Let $\e>0$ be given. By Corollary \ref{corr1}, there exists
$K(\e)=K(\e,{\cal P}(\mu))$ such that \baa
e^{-2K(\e)t}\|u(\cdot,t)\|^2_{H^0}\le
\e\int_0^te^{-2K(\e)s}\|\varphi(\cdot,s)\|^2_{H^0}ds\qquad \forall
t\in(0,T).  \label{aas0}\eaa  Let $\varphi\in X^{-1}_c$. Set \baaa
p_0(\varphi,t)\defi\frac{1}{t}\int_0^t\|\varphi(\cdot,s)\|^2_{H^0}ds,
\quad q(u,t)\defi\|u(\cdot,t)\|_{H^0}.\eaaa It follows that
 \baaa
\sup_{h}\left(\frac{q(u,t)}{tp_0(\varphi,t)}-\frac{1-e^{-2K(\e)t}}{tp_0(\varphi,t)}q(u,t)
\right)\le \e
 \qquad \forall t\in(0,T). \eaaa
 Hence
 \baaa
\sup_{h}\frac{1}{tp_0(\varphi,t)}q(u,t)\le \e+\sup_{h\in
X^0}\frac{1-e^{-2K(\e)t}}{tp_0(\varphi,t)}q(u,t).
 \qquad \forall t\in(0,T). \eaaa
By (\ref{aas0}),
 \baaa q(u,t)\le\e
e^{2K(\e)t} tp_0(\varphi,t) \qquad \forall t\in(0,T). \eaaa Hence
\baaa\sup_{h}\frac{1-e^{-2K(\e)t}}{tp_0(\varphi,t)}q(u,t)\to 0
\quad\hbox{as}\quad t\to 0+\quad \forall  \e>0. \eaaa
\par
If $\varphi\in X^{0}_c$, then $p_0(\varphi,t)\to
\|\varphi(\cdot,0)\|_{H^0}^2$ as $t\to 0+$. It follows that
\baaa\ulim_{t\to 0+}\sup_{\varphi\in
X^{0}_c}\frac{q(u,t)}{t\|\varphi(\cdot,0)\|_{H^0}}\le \e
 \eaaa
for any $\e>0$. Then (\ref{as10}) follows.
 This completes the proof of Theorem \ref{ThAs0}.
 $\Box$
\par {\it Proof of Theorem \ref{ThAs}.}
Let $\e>0$ be given. By Corollary \ref{corr1} again, there exists
$K(\e)=K(\e,{\cal P}(\mu))$ such that \baa
e^{-2K(\e)t}\|u(\cdot,t)\|^2_{H^0}\le
\left(\frac{1}{2}+\e\right)\int_0^te^{-2K(\e)s}(F(\cdot,s),b(\cdot)^{-1}F(\cdot,s))_{H^0}ds
\qquad \forall t\in(0,T),\  \label{aas}\eaa where
$F=(F_1,...,F_n)$, and where $F_i\in X^0$ are such that (\ref{hU})
holds. Let $\varphi\in X^{-1}_c$. Set \baaa
p(F,t)\defi\frac{1}{t}\int_0^t(F(\cdot,s),b(\cdot,s)^{-1}F(\cdot,s))_{H^0}ds,
 \quad q(u,t)\defi\|u(\cdot,t)\|_{H^0}.\eaaa It follows that
 \baaa
\sup_{h}\left(\frac{q(u,t)}{tp(F,t)}-\frac{1-e^{-2K(\e)t}}{tp(F,t)}q(u,t)
\right)\le \left(\frac{1}{2}+\e\right)
 \qquad \forall t\in(0,T). \eaaa
 Hence
 \baaa
\sup_{F:\ \varphi\in X^{-1}_c}\frac{1}{tp(F,t)}q(u,t)\le
\left(\frac{1}{2}+\e\right)+\sup_{h\in
X^0}\frac{1-e^{-2K(\e)t}}{tp(F,t)}q(u,t).
 \qquad \forall t\in(0,T). \eaaa
By (\ref{aas}),
 \baaa q(u,t)\le
e^{2K(\e)t}\left(\frac{1}{2}+\e\right)tp(F,t) \qquad \forall
t\in(0,T). \eaaa Hence \baaa\sup_{F:\ \varphi\in
X^{-1}_c}\frac{1-e^{-2K(\e)t}}{tp(F,t)}q(u,t)\to 0
\quad\hbox{as}\quad t\to 0+\quad \forall  \e>0. \eaaa
\par
If $\varphi\in X^{-1}_c$, then $p(F,t)\to
(F(\cdot,0),b(\cdot,0)^{-1}F(\cdot,0))_{H^0}$ as $t\to 0+$. It
follows that \baaa\ulim_{t\to 0+}\sup_{F:\ \varphi\in
X^{-1}_c}\frac{q(u,t)}{t(F(\cdot,0),b(\cdot,0)^{-1}F(\cdot,0))_{H^0}}\le
\left(\frac{1}{2}+\e\right)
 \eaaa
for any $\e>0$. Then (\ref{as1}) follows.
 This completes the proof of Theorem \ref{ThAs}.
\par
{\it Proof of Theorem \ref{Thn}.}
 Note that $u\in Y^1$ is the solution of the problem (\ref{D})
 if and only if
 $u_K(x,t)\defi e^{-Kt}u(x,t)$ is the solution of the problem  \be \label{D2}\ba \frac{\p
u_K}{\p t} =\A u_K-Ku_K +B_K(u_K)+\varphi_K,\quad \quad t\in(
0,T),\\ u_K(x,0)=0,\quad u_K(x,t)|_{x\in \p D}=0,\ea \ee where
$\varphi_K(x,t)\defi e^{-Kt}\varphi(x,t)$.
 In addition,
 $$ \|u\|_{Y^1}\le
e^{KT}\|u_K\|_{Y_1},\quad \|\varphi_K\|_{X^{-1}}\le
\|\varphi\|_{X^{-1}}. $$ Therefore, the solvability and uniqueness
in $Y^1$ of problem (\ref{D}) follows from existence of $K>0$ such
that problem (\ref{D2}) has an unique solution in $Y^1$. Let us
show that this $K$ can be found.
\par
We introduce  operators $F_K: X^{-1}\to Y^1$ such that
$u=F_K\varphi$ is the solution of problem (\ref{oog3}).
\par
 Let $g\in X^{-1}$ be such that
\baa g=\varphi+B_K(w),\quad  \hbox{where}\quad
w=F_Kg.\label{gw}\eaa  In that case, $u_K\defi F_Kg\in Y^1$ is the
solution of (\ref{D2}). \par Equation (\ref{gw}) can be rewritten
as $g=\varphi+R_K(g)$, or \baa g-R_K(g)=\varphi,\label{gw2}\eaa
where the mapping $R_K:X^{-1}\to X^{-1}$ is defined as
$$R_K(g)=B_K(F_Kg).$$
\par
Let $w=F_Kh$, where $h\in X^{-1}$. By Theorem \ref{ThM}
reformulated as Lemma \ref{lemma<<1}, for  any $\e>0$, $M>0$,
there exists $K(\e,M,{\cal P(\mu)})\ge 0$ and a constant
$C_0=C_0({\cal P(\mu)})$ such that \baa &&\sup_{t\in[0,T]}
\|w(\cdot,t)\|_{H^0}^2+ M\int_0^t \|w(\cdot,s)\|^2_{H^0}ds\le
C_0\|h\|_{X^{-1}}\quad
 \quad \forall h\in
X^{-1}.\label{ineq} \eaa Hence \baaa \|F_K h\|_{X^0}\le
M^{-1}C_0\|h\|_{X^{-1}}.\eaaa Take $M$ and $K$ such that
$\d_*\defi C_*M^{-1}C_0<1$. By (\ref{LipY1}), it follows that
\baa\|R_K(g_1)-R_K(g_2)\|_{X^{-1}}\le
C_0\|F_Kg_1-F_Kg_2\|_{X^{0}}\le \d_*\|g_1-g_2\|_{X^{-1}}.
\label{<1}\eaa By The Contraction Mapping Theorem, it follows that
the equation (\ref{gw2}) has an unique solution $g\in X^{-1}$.
Hence problem (\ref{D2}) has an unique solution $u_K=F_Kg\in Y^1$.
This completes the proof of Theorem \ref{Thn}. $\Box$
\par
{\it Proof of Theorem \ref{ThnY2}}. Let $w=F_Kh$, where $h\in
X^{0}$, and where $F_k$ is the operator  defined in the proof of
Theorem \ref{Thn}. By Lemma 7.1 from Dokuchaev (2008), for any
$\e>0$, $M>0$, there exists $K(\e,M,{\cal P(\mu)})\ge 0$ and a
constant $C_0=C_0({\cal P(\mu)})$ such that \baa
&&\sup_{t\in[0,T]} \|w(\cdot,t)\|_{H^1}^2+ M\int_0^t
\|w(\cdot,s)\|^2_{H^1}ds\le C_0\|h\|_{X^{0}}\quad
 \quad \forall h\in
X^{0}.\label{ineqY2} \eaa The rest  of the proof of Theorem
\ref{ThnY2} repeats the proof of Theorem \ref{Thn} with the
replacement of $Y^1$ for $Y^2$, and $X^{-1}$ for $X^0$, and with
$R_K$ being a mapping $R_K:X^0\to X^0$. $\Box$
\par
{\it Proof of Theorem \ref{ThEx}}. The proof for (i)-(iv)
represents simplified versions  of the proof for (v)-(vi) given
below and will be omitted.
\par
 Let $Q_t\defi \{(y,s)\in Q:\ s\le t\}$. Let us prove
(v). We have that \baaa
 &&|B_K(u_1)(x,t)-B_K(u_2)(x,t)|\\
 &&\le e^{-Kt}\int_{Q_t}
 |\b(e^{Ks}u_1(y,s),x,t,y,s)-\b(e^{Ks}u_2(y,s),x,t,y,s)|dyds\\
 &&\le C_L\int_Q|u_1(y,s)-u_2(y,s)|dy ds\le
C_L\ell_{n+1}(Q)^{1/2} \|u_1(\cdot)-u_2(\cdot)\|_{X^{0}}
 \eaaa for all $u_1(\cdot),u_2(\cdot)\in X^{0}$. Since the domain
 $Q$ is bounded, we have that
 \baaa
 \|B(u_1)-B(u_2)\|_{X^{-1}}\le  \|B(u_1)-B(u_2)\|_{X^0}\le \ell_{n+1}(Q)^{1/2}
  \|B(u_1)-B(u_2)\|_{L_{\infty}(Q)}.
  \eaaa
Hence (\ref{LipY1}) holds.
\par
Further, it follows from the assumptions that $B(0)\in X^0$.
 Hence
 $B(u)\in X^0$ for all $u\in X^0$.
 This completes the proof of statement (v).
\par
Let us prove (vi). By the definition, $B(u)=\nabla\cdot \w B(u)$,
where $\w B:X^{0}\to X^0$ is a mapping similar to the one from
statement (v). Then the proof is similar to the proof of statement
(v).
\par
Let us prove statement (vii). Let us assume that $\w\b =0$.  We
have that \baaa && \|B_K(u_1)-B_K(u_2)\|^2_{X^{-1}}\\&&\le
\int_0^T
e^{-2Kt}\|\b(e^{K\tau(t)}u_1(\cdot,\tau(t)),\cdot,\tau(t))-\b(e^{K\tau(t)}u_2(\cdot,\tau(t)),\cdot,\tau(t))\|_{H^0}^2dt\\
&&+ \int_0^T e^{-2Kt}
\|\w\b(e^{K\tau(t)}u_1(\cdot,\tau(t)),\cdot,\tau(t))-\w\b(e^{K\tau(t)}u_2(\cdot,\tau(t)),\cdot,\tau(t))\|_{H^0}^2dt
\\&&\le 2C_L^2
\int_0^T \|u_1(\cdot,\tau(t)))-u_2(\cdot,\tau(t))\|_{H^0}^2dt\\&&
= 2C_L^2 \int_0^T
\|u_1(\cdot,\tau(t))-u_2(\cdot,\tau(t))\|_{H^0}^2
\left(\frac{d\tau(t)}{dt}\right)^{-1}d\tau(t)
\\&&\le \d_*2C_L^2
\int_{\tau(\t)}^{\tau(T)} \|u_1(\cdot,s)-u_2(\cdot,s)\|_{H^0}^2 ds
\le \d_*2C_L^2 \|u_1-u_2\|_{X^0}^2. \eaaa  By the assumptions, it
follows that  $B(0)\in X^{-1}$. Hence
 $B(u)\in X^{-1}$ for all $u\in X^0$.
\par
Let us prove  statement (viii). We have that $B_K(u)=B(u)$, i.e.,
it is independent from $K$. Further,
 \baaa\label{Bjump2} &&B(u)=\w B(u)+\ww B(u),
 \\&&(\w B(u))(x,t)=\int_{\R^n}\Ind_{\{x+c(x,y,t)\in D\}}
 u(x+c(x,y,t),t)\rho(y,t)dy,\\&&(\ww B(u))(x,t)=
 -u(x,t)\int_{\R^n}\Ind_{\{x+c(x,y,t)\in D\}}\rho(y,t)dy
 \\&&\hphantom{xxxxxxxxxxxxxxxx}
 -\left(\int_{\R^n}\Ind_{\{x+c(x,y,t)\in D\}}c(x,y,t)\rho(y,t)dy\right)^\top\nabla u(x,t).
 \eaaa
  It follows from the assumptions that $\ww B:X^0\to X^{-1}$ is a linear and continuous operator.
  Hence it suffices to prove
 that (\ref{LipY1}) holds for the operator $\w B$.
We have that \baaa (\w B(u))(x,t)=\int_D
 u(z,t)r(x,z,t)dz.
 \eaaa
 Clearly, $\w B(0)=0$. Further, we have that \baaa \|\w B(u_1)-\w B(u_2)\|^2_{X^{-1}}
 \le\|\w B(u_1)-\w B(u_2)\|^2_{X^{0}}=\int_Q \left(\int_D
 (u_1(z,t)-u_2(z,t))r(x,z,t)dz\right)^2 dxdt
\\
\le
\int_Q \left(\int_D
 |u_1(z,t)-u_2(z,t)|^2dz\right) \left(\int_D |r(x,z,t)|^2dz\right) dxdt
\\
 \le
\int_0^Tdt \left(\int_D
 |u_1(z,t)-u_2(z,t)|^2dz\right) \int_Ddx\int_D |r(x,z,t)|^2dz
\\ \le\left(\esssup_{t\in [0,T]} \int_{D\times D}
|r(x,z,t)|^2dxdz\right)\,
   \|u_1-u_2\|_{X^0}^2. \eaaa   This completes the proof of statement (viii)
 and the proof of Theorem \ref{ThEx}. $\Box$
\par
{\it Proof of Theorem \ref{ThSharpM}.}
 Repeat that
$u(x,t)=e^{Kt}u_K(x,t)$, where $u$ is the solution
 of problem  (\ref{4.1})  and $u_K$ is the solution of (\ref{oog3}) for
$h_K(x,t)=e^{-Kt}h(x,t)$. Therefore,  it suffices to find $n$,
$D$, $b,f,\lambda$,  such that\baa &&\forall T>0,c>0, K>0\quad
\exists \varphi\in X^{-1}:\quad  \nonumber\\
&&\|u(\cdot,T)\|_{H^0}^2
\ge\left(\frac{1}{2}-c\right)\int_0^T(F(\cdot,t),b(\cdot,t)^{-1}F(\cdot,t))_{H^0}dt,
\label{prot10}\eaa  where $u$ is the solution
 of problem  (\ref{oog3}) and $F_i\in X^0$ are such as presented in (\ref{hU}), $F=(F_1,...,F_n)$.
\par
 Let us show that (\ref{prot10}) holds for \be
 n=1,\quad D=(-\pi,\pi),\quad b(x,t)\equiv 1,
 \quad f(x,t)\equiv 0,\quad \lambda(x,t)\equiv 0.\label{mu+}\ee
\par
 In this case,
(\ref{oog3}) has the form \baa u'_t=u''_{xx}-Ku+h,\quad
u(x,0)\equiv 0,\quad u(x,t)|_{x\in\p D}=0, \label{eu}\eaa
\par
Let
 \baa
 \g=m^2+K,\quad
\varphi_m(x,t)\defi m\sin(m x) e^{\g t}, \quad
F_m(x,t)\defi-\cos(m x) e^{\g t}, \quad\label{wxi}\eaa where
$m=1,2,3,\ldots$.
 It can be verified immediately that  the
solution of the boundary value problem is \baaa
u(x,t)=m\sin(mx)\int_0^te^{-\g(t-s)+\g s}ds=m\sin(m x)e^{-\g
t}\int_0^te^{2\g s}ds =m\sin(m x)e^{-\g t}\frac{e^{2\g t}-1}{2\g
}. \eaaa Hence \baaa\|u(\cdot,T)\|^2_{H^0}= m^2\|\sin(m
x)\|^2_{H^0}e^{-2\g T}\left(\frac{e^{2\g }-1}{2\g}\right)^2
=m^2\pi e^{-2\g T}\frac{(e^{2\g T}-1)^2}{4\g^2}, \eaaa and \baaa
\int_0^T\|F_m(\cdot,t)\|^2_{H^0}dt=\|\cos(m
x)\|^2_{H^0}\int_0^Te^{2\g t}dt= \pi\frac{e^{2\g T}-1}{2\g}. \eaaa
It follows that \baa
\|u(\cdot,T)\|^2_{H^0}\left(\int_0^T\|F_m(\cdot,t)\|^2_{H^0}dt\right)^{-1}
=\frac{m^2}{2\g}e^{-2\g T}(e^{2\g T} -1)=\frac{m^2}{2\g}(1-e^{-2\g
T})\to \frac{1}{2} \label{1/4}\eaa as $\g\to +\infty$. In
particular, it holds if $K$ is fixed and $m\to +\infty$. It
follows that (\ref{prot1}) holds.
 This completes the proof of Theorem \ref{ThSharpM}. $\Box$
\par
{\it Proof of Theorem \ref{ThSharp}.} Let the parameters be
defined by (\ref{mu+}). Consider a sequence $\{T_i\}$ such that
$T_i\to 0+$ as $i\to +\infty$.  Let $\varphi=\varphi_{m}$ be
defined by (\ref{wxi}) for an increasing sequence of integers
$m=m_i$ such that $m_i>T_i^{-1}$. In that case, $\g T\to +\infty$.
Hence(\ref{1/4}) holds and (\ref{as2}) holds.  $\Box$
\subsection*{Acknowledgment}  This work  was supported by NSERC
grant of Canada 341796-2008 to the author.
\section*{References}
$\hphantom{xx}$
\par
Ammann, H. (2005). Non-local quazi-linear parabolic equations.
{\it Russ. Math. Surv.} {\bf 60} (6), 1021-1033.

Dokuchaev, N.G. (2008). Universal estimate of the  gradient for
parabolic equations. {\it Accepted to Journal of Physics A:
Mathematical and Theoretical}.
\par
Ladyzhenskaia, O.A. (1985). {\it The Boundary Value Problems of
Mathematical Physics}. New York: Springer-Verlag.
\par
Ladyzenskaya, O.A., Solonnikov, V.A.,  and  Ural'ceva, N.N. {\it
Linear and Quasi--Linear Equations of Parabolic Type.} Moscow,
Nauka, 1967. (in Russian) (English translation: Providence, R.I.:
American Mathematical Society, 1968).
\par
Zheng, S. (2004). {\it  Nonlinear evolution equations}.
(Monographs and Surveys in Pure and Applied Mathematics, vol. 133)
Chapman and Hall/CRC, Boca Raton, FL.
\end{document}